\documentclass[11pt]{article}
\usepackage{geometry}
\usepackage{graphicx}
\usepackage{amsmath}
\usepackage{amssymb}
\usepackage{amsthm}
\usepackage{epsfig}
\usepackage{float}
\usepackage{epstopdf}
\usepackage{color}
\usepackage{colordvi}
\usepackage[ruled, vlined]{algorithm2e}
\usepackage{algpseudocode}
\usepackage{subcaption}
\usepackage{siunitx}
\usepackage{comment}
\usepackage{xcolor}
\usepackage{tikz}
\usepackage{pgfplots}
\usepackage{placeins}
\usetikzlibrary{calc}
\usepackage{authblk}
\usepackage[
colorlinks=true,
urlcolor=blue,
linkcolor=green
]{hyperref}
\usepackage[capitalize]{cleveref}
\renewcommand{\Vec}[1]{{\boldsymbol{#1}}}

\newcommand{\bM}{\boldsymbol{M}}
\newcommand{\bK}{\boldsymbol{K}}

\newcommand{\bF}{\boldsymbol F}
\newcommand{\bQ}{\boldsymbol Q}

\newcommand{\bV}{\boldsymbol V}

\newcommand{\bq}{\boldsymbol{q}} 

\newcommand{\bs}[1]{\boldsymbol{#1}}

\newcommand{\hM}{\widehat{\boldsymbol{M}}}
\newcommand{\hK}{\widehat{\boldsymbol{K}}}
\newcommand{\hq}{\widehat{\textbf{q}}}
\newcommand{\hf}{\widehat{\textbf{f}}}

\newcommand{\R}{\mathbb{R}}
\newcommand{\Mat}[1]{{\boldsymbol{#1}}}

\usepackage{todonotes}

\usepackage[normalem]{ulem}

\graphicspath{{Plots/}}

\title{Inference of Substructured Reduced-Order Models for Dynamic Contact from Contact-free Simulations}

\author[*]{Diana Manvelyan-Stroot} 
\author[$\dagger$]{Yevgeniya Filanova}   
\author[$\dagger\dagger$]{Igor Pontes Duff}
\author[$\ddagger$]{Peter Benner}
\author[$\ddagger \ddagger$]{Utz Wever}

\affil[*]{\small Siemens AG,  Munich, diana.manvelyan@siemens.com}
\affil[$\dagger$]{ Otto-von-Guericke University Magdeburg, Max Planck Institute for Dynamics of Complex Technical Systems, filanova@mpi-magdeburg.mpg.de}
\affil[$\ddagger$]{Max Planck Institute for Dynamics of Complex Technical Systems, Otto-von-Guericke University Magdeburg,  benner@mpi-magdeburg.mpg.de}
\affil[$\dagger\dagger$]{ Max Planck Institute for Dynamics of Complex Technical Systems,  pontes@mpi-magdeburg.mpg.de}
\affil[$\ddagger \ddagger$]{utz.wever@t-online.de}

\date{}

\pgfplotsset{compat=1.18}

\begin{document}
\maketitle
\begin{abstract}
In this paper, we propose an operator-inference-based reduction approach for contact problems, leveraging snapshots from simulations without active contact. Contact problems are solved using adjoint methods, by switching to the dual system, where the corresponding Lagrange multipliers represent the contact pressure. The Craig-Bampton-like substructuring method is incorporated into the inference process to provide the reduced system matrices and the coupling of the contact and interior nodes. The maximum possible set of contact nodes must be known a priori.
Characteristic properties of the inferred matrices, such as symmetry and positive definiteness, are enforced by appending additional constraints to the underlying least-squares problem. The resulting dual system, which forms a linear complementarity problem, is well-defined and can be effectively solved using methods such as Lemke’s algorithm. The performance of the proposed method is validated on three-dimensional finite element models.

\end{abstract}



\newpage
\section{Introduction}\label{sec:1}
Nowadays, digital twins are valuable commodities in the technology landscape, offering enormous benefits such as online prediction and optimization of complex time-dependent industrial systems. Despite their potential, fully exploiting their capabilities requires the combination and integration of multiple technologies across different disciplines. This remains a significant challenge and is a focus of ongoing research \cite{van2022executable}.

One of the benefits of current research is that the individual disciplines, that are needed to construct the digital twins, are already well-established \cite{hartmann2023real}. One of these important disciplines is model order reduction \cite{benner2005, hartmann202012}, which addresses the challenge of providing small-scale models that are real-time capable while maintaining the essential physical characteristics relevant to specific applications. In other words, scalable approaches are required for enabling digital twins.
Extensive research has been done on model order reduction, see \cite{Baur2014,lohmann2004,Farhat2014,bennerwillcox} including data-driven methods \cite{radermacher2016pod, peherstorfer2016data, swischuk2019projection, filanova2023operator} and realization-based methods such as balanced truncation \cite{benner2016frequency, besselink2014model, kurschner2018balanced}, Krylov subspace methods \cite{ilyas2016krylov, breiten2010krylov, benner2011recycling} and modal reduction \cite{Baur2014}.

In this work, we address digital twins for applications in structural mechanics, in particular for dynamic contact, which are used in the context of operational support simulation facilitating improved maintenance scheduling, life expectancy estimations, advanced fault detection, and operational control and optimization \cite{heinrich2019online, martynov2019polynomial}. For further background literature on contact mechanics including discussions on variational inequalities and classical numerical methods such as the penalty method or the augmented Lagrangian formulation, we refer to \cite{Fran75,KiOd88,wriggers2004computational}. Although these methods lead to unconstrained minimization problems, they often prevent the accurate determination of the dual variable and the preservation of the contact shape.

Our research is focused on the development of reduced models that can provide accurate solutions for displacements and contact pressures. In this context, the Lagrange multiplier method is particularly well-suited, wherein the Lagrange multipliers represent the discrete contact pressure. One of the earliest reduction methods in relation to the Lagrange multiplies is a snapshot-based method that uses the singular value decomposition and the non-negative matrix compression to reduce all quantities of interest for a linear dynamic contact problem \cite{Charbel2017}. Another way to reduce the Lagrange multipliers is proposed in many hyper-reduction methods based on a reduced integration domain \cite{fauque2018hybrid, le2022condition}.

For dynamic linear-elastic contact problems model order reduction, approaches based on Krylov subspace methods and adjoint methods, were introduced in \cite{manvelyan2021efficient, manvelyan2022physics}. These approaches are based on the primal-dual decoupling method and Craig-Bampton substructuring, which leads to the reduction only of the primal displacements while maintaining complete information within the discretized contact area. This is particularly advantageous when the contact area is minor compared to the overall structure. After the reduced primal system is defined, we switch over from the space of reduced displacements to the space of the dual Lagrange multipliers, i.e., to the adjoint problem of the original variational problem. This leads to a linear complementarity problem (LCP) \cite{CotPetal2009,Li2016}. The complexity of the LCP scales with the number of constraints, which is usually small compared to the number of displacements and thus, may be solved even in real-time. The resulting update of the displacements is again performed in a reduced space. The degree of reduction is optimal in the sense that the shape of the contact interface is preserved.

The above mentioned approaches require the extraction of the system matrices of the semi-discretized dynamic equations, which can be a challenge for modern finite element simulation software due to the inaccessibility of the implemented solution algorithms. Therefore, we focus on non-intrusive reduction approaches that can provide the reduced models using the simulation data only. Generating the snapshots is a crucial part that influences the efficiency of the reduction since it requires the time-consuming integration of the high-fidelity model. Therefore, we propose to use less cumbersome contact-free simulations for reduced order modeling, i.e., using only displacement and force snapshots from the simulations without active contact. Avoiding contact simulations for data generation allows us to save computational resources and extend the range of suitable software options.
In this work, we first learn the primal system without contact, preserving its second-order structure via operator inference procedure as described in \cite{morPehW16}.
Additionally, a substructuring in the fashion of Craig-Bampton \cite{craigbampton} is incorporated into the learning procedure of the primal system matrices. Here, special attention is devoted to restoring the coupling information between the interior and boundary displacement nodes. The latter group of nodes represents the maximum possible set of contact nodes and has to be defined after generating the contact-free simulation data. For that, we suggest three different approaches: a specific least-square approach, one in full and one in reduced space, and a data-based method by means of successive static simulations \cite{subAllet2020, delhez2021reduced}. The choice between these three approaches is a matter of trade-off between accuracy and efficiency and will be discussed in more detail throughout this paper. Similar to \cite{manvelyan2021efficient}, here we utilize the primal-dual decoupling method and switch to the adjoint system resulting in a LCP. As it turns out, physical properties of the primal system matrices such as symmetry and positive definiteness are necessary for solvability of the obtained LCP system. Hence, the extension of the
operator inference procedure to a constrained convex optimization problem \cite{filanova2023operator} is used. Finally, the obtained reduced LCP is solved by an appropriate method, e.g., Lemke's algorithm, see \cite{manvelyan2021efficient}.

The main novelty of our method can be summarized as follows. The operator inference procedure with LMI constraints is tailored with a substructuring technique for contact treatment. The resulting ROM is solved with adjoint methods allowing accurate computation of both the displacement and the contact pressure. 

The paper is organized as follows: In \cref{sec:2} we briefly recall the underlying equations of dynamic contact in linear elasticity and the corresponding solution based on adjoint methods. The \cref{sec:3} presents the current reduction methods, used as main tools for the novel reduced order modeling, presented in \cref{sec:4}. The performance of the approach is demonstrated by means of two carefully chosen examples in \cref{sec:5}.


\section{Dynamic contact}\label{sec:2}
This section addresses modeling aspects of large-scale dynamic node-to-node contact models as well as their solution procedures involving Lagrange multiplier formulation and adjoint methods \cite{manvelyan2021efficient}. Our ultimate goal is to introduce a data-driven reduced-order model for the semi-discretized equations that describe 
the nodal displacement variable, denoted by $\Vec{q}$ within a much smaller space while explicitly preserving
the constraints and the dimension of the Lagrange multipliers.  

\subsection{Dynamic Contact Model}
First, we summarize the governing equations for dynamic contact in linear elasticity, presenting both the strong and weak formulations of the obstacle problem. Following the weak formulation, we provide a brief explanation of the node-to-node contact condition formulation.

\subsubsection{Obstacle Problem}

Throughout this work, we consider frictionless, adhesive-free normal contact in combination with small deformation theory and a linear-elastic material. The dynamic contact problem is expressed in terms 
of the displacement field $\Vec{u}(\Vec{x},t)\in\R^d$ and the contact pressure 
$p(\Vec{x},t)\in \R$ where $x \in \Omega \subset \R^d, d=2$ or $d=3$, is the spatial variable 
and $t \in [t_0,T]$ stands for the temporal variable.
The boundary of the elastic body is decomposed into
$\partial\Omega = \Gamma_D\cup\Gamma_N\cup\Gamma_C$ where the  
latter represents the contact interface. Using the outward normal vector $\Vec{n}(\Vec{x})$ on
the contact interface, the distance between the body and a given surface is  described by the scalar gap function 
\begin{equation}\label{eq:defgap}
g(\Vec{x}, \Vec{u}) = \big(\Vec{\xi}(\Vec{x}) - (\Vec{x} + \Vec{u}(\Vec{x},t)\big)^T \Vec{n}(\Vec{x}), \quad \Vec{x} \in \Gamma_C.
\end{equation}
Here, the point $\Vec{\xi}(\Vec{x})$ on the obstacle's surface is obtained by projection of $\Vec{x}$ in the outward normal direction.
 The contact problem in strong form is then given by the dynamic equations
\begin{equation}\label{elast_cont} 
\rho\Mat{u}_{tt}-\mbox{div} \Mat{\sigma}(\Vec{u})= \Vec{F} \quad \mbox{in } \Omega \times [t_0,T]
\end{equation}
subject to standard boundary conditions
\begin{equation}\label{bcDN}
\Mat{\sigma}(\Vec{u})\cdot\Vec{n}=\Vec{\tau} \quad \mbox{on } \Gamma_N \times [t_0,T], \qquad
\Vec{u}=0  \quad \mbox{on } \Gamma_D \times [t_0,T]
\end{equation}
and furthermore subject to the contact conditions 
\begin{equation}\label{bcContact} 
g \geq 0, \quad
p \geq 0, \quad g \cdot p = 0  \quad \mbox{on } \Gamma_C \times [t_0,T].
\end{equation}
As initial data, 
\begin{equation}
\Mat{u}(\cdot,t_0) = \Mat{u}_0 \quad \mbox{and } \,
\Mat{u}_t(\cdot,t_0) = \Mat{v}_0
\end{equation}
are prescribed with given functions $\Vec{u}_0$ and $\Vec{v}_0$, respectively. 

In (\ref{elast_cont}) and (\ref{bcDN}), the constant $\rho$ denotes the mass density,
$\Vec{F}(\Vec{x},t) \in \R^d$ the volume force,
$\Vec{\tau}(\Vec{x},t) \in \R^d$ the surface traction, 
and 
$\Mat{\sigma}(\Vec{u}) \in R^{d\times d}$ the stress tensor given by
\begin{equation}\label{stress}
  \Mat{\sigma}(\Vec{u}) = \frac{E}{1+\nu} \Mat{e}(\Vec{u}) + \frac{\nu E}{(1+\nu)(1-2 \nu)} \mbox{trace}(\Mat{e}(\Vec{u})) \Mat{I}
\end{equation}
with Young's modulus $E \geq 0$, Poisson's ratio $-1\leq \nu \leq 0.5$ and linearized
strain tensor 
\begin{equation} \label{strain}
  \Mat{e}(\Vec{u})=\frac{1}{2}(\nabla \Vec{u}+\nabla \Vec{u}^T) \in \R^{d\times d}. 
\end{equation}

\subsubsection{Weak Formulation}

In order to pass to a weak formulation of the obstacle problem, we introduce
the function spaces
\begin{eqnarray}
\mathcal{V} &:= & \left\{ \Vec{v} \in H^1(\Omega)^d :
                           \Vec{v} = 0 \, \mbox{on } \Gamma_D \right\}, \label{defV} \\
\mathcal{L} &:= & \left\{ \mu \in H^{1/2}(\Gamma_C)^\prime : 
                         \int_{\Gamma_C} \mu w \, ds \geq 0 \,\, \forall
                          w \in H^{1/2}(\Gamma_C), w \geq 0 \right\}
                          \label{defL}
\end{eqnarray}
and the abstract notation
\begin{equation}\label{defANot}
\langle \rho \Vec{u}_{tt}, \Vec{v} \rangle := \int_\Omega \Vec{v}^T \rho \Vec{u}_{tt} dx, \,\,
a(\Vec{u}, \Vec{v}) := \int_\Omega\Mat{\sigma}(\Vec{u}):\Mat{e}(\Vec{v})dx,
\,\,
\langle \ell , \Vec{v} \rangle := \int_\Omega \Vec{v}^T \Mat{F}dx + \int_{\Gamma_N} \Vec{v}^T\Vec{\tau}ds.
\end{equation}
The non-penetration condition in weak form with test function $\mu \in {\cal L}$ is recast as 
\begin{equation}\label{contactweak}
 0 \leq \int_{\Gamma_C} g \mu ds = 
        \int_{\Gamma_C} g (\mu-p) ds  
   = -\int_{\Gamma_C} \Vec{u}^T\Vec{n} (\mu-p) ds
          + \int_{\Gamma_C} (\Vec{\xi}-\Vec{x})^T\Vec{n} (\mu-p) ds.
\end{equation}
The last two integrals give rise to the definitions of the bilinear form
\begin{equation}\label{defBilb}    b_c(\Vec{v},p) :=  - \int_{\Gamma_C}  \Vec{v}^T \Vec{n} \cdot p \,  ds
\end{equation}
on ${\cal V} \times {\cal L}$
and the linear form on $\cal L$ 
\begin{equation}\label{defLFm}
\langle {m}, p \rangle := \int_{\Gamma_C} (\Vec{\xi}-\Vec{x})^T\Vec{n} p \, ds.
\end{equation}
Using these definitions, the weak form of the dynamic contact problem is stated as follows: For each $ t \in [t_0, T] $ find the displacement field $ \Vec{u}(\cdot,t) \in \mathcal{V} $ and the contact pressure $ p(\cdot,t) \in \mathcal{L} $ such that  
\begin{equation}\label{elast_weak}
\begin{aligned} 
\langle \rho \Vec{u}_{tt}, \Vec{v} \rangle + a(\Vec{u}, \Vec{v}) & = \langle \ell,\Mat{v} \rangle  + b_c(\Vec{v}, p) \quad &&\mbox{for all } \Vec{v} \in \mathcal{V},\\
b_c(\Vec{u}, \mu-p) +  \langle {m}, \mu-p \rangle & \geq 0 \quad &&\mbox{for all }  \mu \in \mathcal{L}. 
\end{aligned}
\end{equation}
The system \eqref{elast_weak} is discretized with respect to the spatial variable by applying the standard Galerkin projection $\Vec{u}(\Vec{x},t) \doteq \sum_{i=1}^N \Vec{\phi}_i(\Vec{x}) q_i(t)
= \Mat{\Phi}(\Vec{x}) \Vec{q}(t)$ with basis functions $\Vec{\phi}_i$ and
nodal variables $\Vec{q} = (q_1, \ldots, q_N)$ to the displacement field.
For a contact problem with $k$ bodies that fits into the framework described so far,
the mass matrix $\Vec{M}$ will consist of $k$ blocks on the diagonal that  stem from
the discretizations of the individual bodies, and the same block-diagonal structure applies to the stiffness matrix $\Vec{K}$. 
If the contact condition is simply expressed in terms of the distance between a node
$\Vec{x}_j \in \Gamma_C, i=1,\ldots,m$, and the obstacle, the integral over $\Gamma_C$ 
in the weak dynamic equation
is replaced by a sum
\begin{equation}\label{discretGammaC}
     \int_{\Gamma_C}  \Vec{v}^T \Vec{n} \cdot p \,  ds \doteq
                     \sum_{j=1}^{m}  A_j \Vec{v}(\Vec{x}_j)^T \Vec{n}(\Vec{x}_j) \cdot p(\Vec{x}_j,t), \quad A_j: \mbox{area around node }\Vec{x}_j.
\end{equation}
Setting $\Vec{\lambda}(t) := (p(\Vec{x}_1,t), \ldots, p(\Vec{x}_{m},t))$ as discrete pressure variable (the geometrical scaling factors are ignored here),
follow in the usual way, with the matrix $\Mat{C} \in \R^{m \times N}$ being an indicator matrix for the nodes in the contact interface and $\Vec{b}\in\R^m$ the offset of the contact. Note that only a few displacements are involved for the contact condition and hence the matrix $\Mat{C}$ exhibits a sparse structure. We do not dive further into the details of various contact models and discretization schemes and refer instead to the standard references \cite{KiOd88,wriggers2004computational}.

The Lagrange multiplier method is used to enforce the non-penetration condition and the finite element method on matching meshes, i.e., node-to-node contact condition is applied for the discretization in space.
Other Approaches such as the augmented Lagrange method or the Nitsche method 
\cite{wriggers2004computational} lead to modifications 
that will not be considered here.

\subsection{Contact Algorithm using Adjoint Methods}
\label{subsec:primaldual}
Popular methods to solve the mechanical contact problem are the penalty or the augmented Lagrange method \cite{Gill81,wriggers2004computational}.
In this work we advocate the Lagrange multiplier method and the primal-dual decoupling method from \cite{manvelyan2021efficient}, which  leads to a LCP.

In the stationary case, the following minimization problem is recast
\begin{equation}\label{eq:var}
\min_{\Vec{q}\in\R^N,\Vec{\lambda}\in\R_+^m} \frac{1}{2}\Vec{q}^T \Mat{K} \Vec{q} - 
\Vec{q}^T \Vec{f} - \Vec{\lambda}^T (\Mat{C}\Vec{q}+\Vec{b}).
\end{equation}

In the transient case, the minimization problem \eqref{eq:var} has to be extended accordingly and the corresponding KKT conditions lead to the dynamic equations
\begin{eqnarray}\label{eq:KKT_transient1}
&&\Mat{M}\ddot{\Vec{q}}(t) + \Mat{K} \Vec{q}(t)  = \Vec{f}(t)  + \Mat{C}^T\Vec{\lambda}(t), \\ \label{eq:KKT_transient2}
&&\Mat{C} \Vec{q}(t)  + \Vec{b}\geq 0, \quad \Vec{\lambda}(t) \geq 0, 
	\quad \Vec{\lambda}(t)^{T} (\Mat{C}\Vec{q}(t) +\Vec{b})=0
\end{eqnarray}

As a straightforward and unconditionally stable method, we apply the implicit 
Euler scheme, i.e., the first-order backward differentiation formula. 
For a given time stepsize $h$, the second-order derivative is replaced by the finite difference
\begin{equation}\label{eq:findiff}
\ddot{\Vec{q}}(t+h) \approx \frac{1}{h^2}\big(\Vec{q}(t+h) - 2\Vec{q}(t) + \Vec{q}(t-h)\big),
\end{equation}

Inserting \eqref{eq:findiff} into \eqref{eq:KKT_transient1} leads to
\begin{equation}\label{eq:system_disc}
\Mat{M}\left(\Vec{q}(t+h) -2\Vec{q}(t) + \Vec{q}(t-h)\right) + h^2\Mat{K} \Vec{q}(t+h) = 
h^2\Vec{f}(t+h) + h^2\Mat{C}^T\Vec{\lambda}(t+h)
\end{equation}
Using the fact that the mass and stiffness matrices are positive-definite and assuming that the solution at the previous time steps are known, the dynamic equation \eqref{eq:system_disc} may be resolved for $\Vec{q}(t+h)$ as follows
\begin{equation}\label{eq:system_disc_resolved}
\Vec{q}(t+h) = (\Mat{M}+h^2\Mat{K})^{-1}(h^2\Vec{f}(t+h) + h^2\Mat{C}^T\Vec{\lambda}(t+h) +2\Mat{M}\Vec{q}(t) - \Mat{M}\Vec{q}(t-h)).
\end{equation}
Inserting (\ref{eq:system_disc_resolved}) into the constraints (\ref{eq:KKT_transient2}) leads to an LCP
\begin{equation} \label{eq:KKTAB}
	\begin{array}{rcc}
	\Vec{B} + \Mat{A} \Vec{\lambda} & \geq & 0, \\
	\Vec{\lambda} & \geq & 0, \\
	\Vec{\lambda}^T (\Vec{B} + \Mat{A} \Vec{\lambda}) & = & 0.
	\end{array}
\end{equation}
where
\begin{eqnarray} \label{eq:AB_LCP}
\Mat{A} &:=& h^2\Mat{C}(\Mat{M}+h^2\Mat{K})^{-1}\Mat{C}^T,\\
\Mat{B} &:=& \Mat{C}(\Mat{M}+h^2\Mat{K})^{-1}(h^2\Vec{f}(t+h) + 
        2\Mat{M}\Vec{q}(t) - \Mat{M}\Vec{q}(t-h)) + \Vec{b}
\end{eqnarray}
The LCP problem (\ref{eq:KKTAB})  may be solved by standard methods from
 constrained optimization, see below for more details. Note that for a constant time step $h$, the LCP matrix $\Mat{A}$ is independent of time and thus needs to be computed only once. The LCP vector $\Vec{B}$ has to be updated in each time step before starting the LCP solving procedure.

We remark that the two-step time discretization requires initial values
$\Vec{q}(t_0)$ and $\Vec{q}(t_0+h)$ to start. E.g., if $\Vec{q}_0$ and
$\dot{\Vec{q}}_0$ as initial displacement and velocity are given, one can
compute $\Vec{q}(t_0+h) = \Vec{q}_0 + h \dot{\Vec{q}}_0$ by an explicit Euler step and then continue with the two-step formula (\ref{eq:system_disc}).
The Lagrange multiplier $\Vec{\lambda}$, on the other hand, does not require 
an initial value, and it is computed in each time step as an implicitly given function
of $\Vec{q}$. In the terminology of differential-algebraic equations, this means that 
$\Vec{q}$ stands for the differential variables while the algebraic variables 
$\Vec{\lambda}$ possess no memory.

\section{Applied Reduction Methods} \label{sec:3}
Finite element models often have a large number of degrees of freedom (DOFs), making the LCP solution \eqref{eq:KKTAB} very computationally expensive. Hence, it is advantageous to reduce the problem size by approximating the state vector of the high-fidelity model $\Vec{q}$ in a lower-dimensional subspace by the reduced state vector $\widehat{\Vec{q}}$:
\begin{equation}\label{eq:red_q}
    \Vec{q} \approx \Mat{V} \widehat{\Vec{q}}.
\end{equation}
Many intrusive reduced models are generated by the projection of the system matrices onto the subspace spanned by $\Mat{V}$ \cite{morSaaSW19}. In this work, we assume that there is no access to the system operators. Therefore, we develop a procedure to derive the reduced system matrices only from the simulation data.

The novel non-intrusive substructure-based reduction framework is mainly based on the intrusive Hurty-Craig-Bampton reduction, introduced in \cref{subsec:craig-bampton}, and the force-informed operator inference method, briefly presented in \cref{subsec:fi_opinf}.

\subsection{Intrusive Hurty-Craig-Bampton Reduction} \label{subsec:craig-bampton}
In this subsection we describe the reduction technique, proposed by Hurty \cite{Hur1960}, and later by Craig and Bampton \cite{Craig1968}. This method is typically used in the context of substructuring problems, when a structure is decomposed into subcomponents that are analyzed separately \cite{subAllet2020}. In order to re-assemble the components, the nodes are partitioned into boundary and interior groups, with the displacement vectors $\Vec{q}_B$ and $\Vec{q}_I$ respectively. 

The aforementioned technique is also beneficial for contact mechanics problems, considered in this work. Similarly, we distinguish between the boundary DOFs for the nodes in the contact zone, and the remaining interior DOFs, such that the low-dimensional approximation \eqref{eq:red_q} is:
%
%
%
%
%
\begin{equation} \label{mat:red_q_partitioned}
    \Vec{q} = \begin{pmatrix}
    \Vec{q}_B \\
    \Vec{q}_I
    \end{pmatrix} \approx \Mat{V} \begin{pmatrix}
    \Vec{q}_B \\
    \widehat{\Vec{q}}_I
    \end{pmatrix}.
\end{equation}
 Note that the $\Vec{q}_B$ are preserved in their original dimension, while the reduction acts only on the interior subpart, which ensures consistency for the contact pressure distribution of the full- and reduced-order models.

In order to reveal the most important aspects of Hurty-Craig-Bampton reduction, we consider the partitioned semi-discretized system without contact conditions:
\begin{align} \label{contactfreeODE}
\begin{pmatrix}
\Mat{M}_{BB} & \Mat{M}_{BI}\\
\Mat{M}_{IB} & \Mat{M}_{II}
\end{pmatrix} \begin{pmatrix}
    \ddot{\Vec{q}}_B (t) \\
    \ddot{\Vec{q}}_I (t)
    \end{pmatrix} + \begin{pmatrix}
\Mat{K}_{BB} & \Mat{K}_{BI}\\
\Mat{K}_{IB} & \Mat{K}_{II}
\end{pmatrix} \begin{pmatrix}
    \Vec{q}_B (t) \\
    \Vec{q}_I (t)
    \end{pmatrix} = \begin{pmatrix}
\Mat{f}_B (t)\\
\Mat{f}_I (t)
\end{pmatrix},
\end{align}
where the displacement and force vectors are $ \Vec{q} \in \mathbb{R}^N $, $\Vec{f} \in \mathbb{R}^{N}$, the mass, stiffness matrices are $\Mat{M}, \Mat{K} \in \mathbb{R}^{N \times N}$. We assume that the interior displacements can be represented as a sum of two contributions: static and dynamic. The first one results from the static elimination of the interior DOFs by the boundary DOFs, using the coupling term, calculated as:

\begin{equation} \label{mat:coupling}
    \Mat{\Phi}_{IB} = -\Mat{K}_{II}^{-1}\Mat{K}_{IB}.
\end{equation}
The dynamic contribution characterizes the dynamic behavior of the interior subpart and is represented by reduced interior coordinates $\widehat{\bq}_I$ in a reduced basis $\Mat{V}_{I}$, consisting of interior dynamic modes. The static and dynamic components yield the final expression for the full-order interior displacements:

\begin{equation} \label{eq:red_q_relation}
    \bq_{I} = \bs{\Phi}_{IB} \bq_B + \Mat{V}_{I} \widehat{\bq}_I.
\end{equation} The interior basis $\Mat{V}_{I}$ is usually obtained by fixing the interface DOFs and solving the generalized eigenvalue problem for the interior subpart $\Mat{M}_{II}$ and $\Mat{K}_{II}$, such that the first $r$ resulting eigenmodes form the basis.
Exploiting the relations \eqref{eq:red_q_relation} and \eqref{mat:red_q_partitioned} for the reduced and full interior state vector, the global reduction matrix is constructed as

\begin{equation}\label{mat:V}
\Mat{V} =
\begin{pmatrix}
\Mat{Id}_I &  0 \\
\Mat{\Phi}_{IB} & \Mat{V}_{I}.
\end{pmatrix},
\end{equation}
where $\Mat{Id}_I \in \mathbb{R}^{n_B \times n_B} $ is the identity matrix. The ROM operators are then obtained by projecting the system matrices onto the low-dimensional subspace.
%
%
Thus, we get a substructure-based reduced surrogate model ($\widehat{\bM}, \widehat{\bK}$) that can be efficiently applied for contact mechanics problems. 
%
However, the application is restricted to the case, when the system matrices are accessible in their partitioned form.

\subsection{Force-informed operator inference} \label{subsec:fi_opinf}

In contrast to the above mentioned intrusive methodology, the operator inference method can be applied, when access to the system operators is not possible. The method was originally introduced for first-order ODE systems in \cite{morPehW16}. Since then numerous modifications, suited for different problem types, were proposed, see, e.g., in \cite{morKraetal2024, morSha24}. In this work, we use the extension of the operator inference method for second-order ODE systems, namely the force-informed approach, proposed in \cite{filanova2023operator}, which is briefly described below.

The essence of the method is to fit the reduced system operators to the low-dimensional approximation of the simulation data by solving a least-squares optimization problem. First of all, we assume to have access to the simulation data, i.e., the displacements and acceleration resulting from integration of \eqref{contactfreeODE} under reasonable external loading conditions. The integration is performed at the predefined time steps $t_j \in [t_0, t_k]$.
Moreover, the external force at each time point is assumed to be completely known. The displacements $\textbf{q}(t)$, and forces $\textbf{f}(t)$ of the system are extracted into the snapshot matrices: 

\begin{equation}
\label{mat:snap}
\textbf{Q} =
\begin{pmatrix}
| &  \dots & | \\
\mathbf{q}(t_1) & \dots & \mathbf{q}(t_k) \\
| & \dots & |
\end{pmatrix} \;
\in \mathbb{R}^{N \times k}, \quad \textbf{F} =
\begin{pmatrix}
| &  \dots & | \\
\mathbf{f}(t_1) & \dots & \mathbf{f}(t_k) \\
| & \dots & |
\end{pmatrix}
\in \mathbb{R}^{N \times k} .
\end{equation} The second derivatives for the displacements $\ddot{\textbf{q}}(t)$ can be found by numerical differentiation, or obtained directly from the simulation if the solver allows it (e.g., second-order integration methods such as Newmark-$\beta$, Generalized-$\alpha$, etc. \cite[Chapter~7]{Ger2015}). The snapshot matrix of the derivatives $\ddot{\bQ}$ can then be constructed analogous to \cref{mat:snap}. 

Note, that the available data \cref{mat:snap} are high-dimensional, therefore the projection onto a low-dimensional subspace is used for the ROM identification. Usually, the POD subspace of the snapshot matrix $\Mat{Q}$ with the basis $\textbf{V}_q$ is used, which leads to the reduced snapshot data:

\begin{align}
    \widehat{\textbf{Q}} = \textbf{V}_q^T \textbf{Q}, & \quad \ddot{\widehat{\textbf{Q}}} = \textbf{V}_q^T \ddot{\textbf{Q}}, \quad  \widehat{\textbf{F}} = \textbf{V}_q^T \textbf{F}. 
\end{align}

 Finally, the ROM operators are the solution of the following constrained optimization problem:

 \begin{equation} \label{eq:finf_opt}
    \min_{\hM \succ 0, \hK \succ 0} \| \hM \ddot{\widehat{\textbf{Q}}}  + \hK \widehat{\textbf{Q}} - \widehat{\textbf{F}} \| _{_F}^2 ,
 \end{equation} which can be efficiently solved using semidefinite programming tools, see \cite{diamond2016cvxpy, Boy94}. The resulting operators $\hM$ and $\hK$ are SPD due to the enforced constraints, which also guarantees the preservation of the stability properties of the original system \cite{morSalEL06}.

\section{Non-intrusive substructuring-based reduction for dynamic contact} \label{sec:4}
Reducing the system without knowing the system representation is particularly beneficial for contact problems, because extracting the system matrices, especially in their "substructured" form, is usually very cumbersome. Based on the methodologies described in \cref{sec:3}, we propose the non-intrusive substructuring method to construct reduced-order models, that enable the approximation of the system response under the contact conditions. 

\subsection{Basic assumptions} \label{subsec:basic_assump}
The problem, considered in our work, is complex and has many different aspects. For the sake of clarity, we summarize the assumptions that were made to develop the procedure, described in the next subsection.

The main assumptions are:

\begin{itemize}
    \item a linear elasticity problem without obstacle \eqref{contactfreeODE} is solved using the implicit Euler integration scheme;
    \item the actual system representation $(\Mat{M}, \Mat{K})$ is unknown;
    \item the simulation solution $\textbf{q} (t)$ is known for predefined time points $t_j \in [t_0, t_k]$;
    \item the external force  vector $\textbf{f} (t)$ is accessible at each time point $t_j$;
    \item the anticipated contact region is predefined, i.e., the indices of the contact nodes are available;
    \item the displacement and force vector are permuted, independent of the simulation software, such that
\begin{equation}\label{eq:permuted_qf}
    \textbf{q} (t) = \begin{pmatrix}
        \textbf{q}_B (t) \\ \textbf{q}_I (t)
    \end{pmatrix},  \quad \textbf{f} (t) = \begin{pmatrix}
        \textbf{f}_B (t) \\ \textbf{f}_I (t)
    \end{pmatrix},
\end{equation} where $\textbf{q}_I$, $\textbf{f}_I \in \mathbb{R}^{n_I}$, and $\textbf{q}_B$, $\textbf{f}_B \in \mathbb{R}^{n_I}$.
    \item the size of the contact zone is small in comparison to the interior part of the structure, i.e., $n_B \ll n_I$.
\end{itemize}


Given the above mentioned assumptions, we aim to obtain the reduced-order model that satisfies

\begin{equation} \label{eq:red_system}
\hM \ddot{\hq} (t) + \hK \hq (t) = \hf (t),    
\end{equation} 
where $\hM$, $\hK \in \mathbb{R}^{(r+n_B) \times (r+n_B)}$ represent the symmetric and positive reduced system operators. The reduced state vector $\hq = \left( \bq_B^T, \; \hq_I^T \right)^T$ includes boundary displacement from the full-order model and the reduced interior displacements $\hq_I \in \mathbb{R}^{r}$. Analogous to the intrusive Hurty-Craig-Bampton, the system \eqref{eq:red_system} stems from the projection of \eqref{contactfreeODE} onto a subspace spanned by $\bV$ presented in \eqref{mat:V}. The ROM force vector is calculated as $\widehat{\textbf{f}} (t) = \bV^T \textbf{f} (t) $.
%


\subsection{Inferring the primal system matrices} \label{subsec:main_procedure}
In the following subsection we will describe the novel methodology to identify the reduced system operators according to the assumptions, specified in the previous subsection. The proposed method consists of three main steps.  
 
 The first step is to arrange the simulation data, obtained as discussed in \cref{subsec:basic_assump}, into the snapshot matrices as follows:

\begin{align}
\label{mat:snap_Q}
\textbf{Q}_B =
\begin{pmatrix}
| &  \dots & | \\
\mathbf{q}_B(t_1) & \dots & \mathbf{q}_B(t_k) \\
| & \dots & |
\end{pmatrix}, \quad \textbf{Q}_I =
\begin{pmatrix}
| &  \dots & | \\
\mathbf{q}_I(t_1) & \dots & \mathbf{q}_I(t_k) \\
| & \dots & |
\end{pmatrix},
\end{align}
 where $ \bQ_B
\in \mathbb{R}^{n_B \times n_k}$, $\bQ_I \in \mathbb{R}^{n_I \times n_k}$. Additionally, the force snapshot matrices are defined analogously:
\begin{align}
\label{mat:snap_F}
\textbf{F}_B =
\begin{pmatrix}
| &  \dots & | \\
\mathbf{f}_B(t_1) & \dots & \mathbf{f}_B(t_k) \\
| & \dots & |
\end{pmatrix}, \quad \textbf{F}_I =
\begin{pmatrix}
| &  \dots & | \\
\mathbf{f}_I(t_1) & \dots & \mathbf{f}_I(t_k) \\
| & \dots & |
\end{pmatrix},
\end{align}
 where $ \bF_B
\in \mathbb{R}^{n_B \times n_k}$, $\bF_I \in \mathbb{R}^{n_I \times n_k}$.

The second step is the calculation of the reduced basis $\bV$, as shown in \eqref{mat:V}, that consists of the coupling term $\Mat{\Phi}_{IB}$ and the interior reduced basis $\Mat{V}_{I}$. The calculation of these two components involves access to the system operators, which is unsuitable for this work. We propose to determine  $\Mat{\Phi}_{IB}$ and $\Mat{V}_{I}$ in a non-intrusive way, as discussed in \cref{subsubsec:interior_basis} and \cref{subsubsec:coupling_mat} respectively.

The final step is the assembling of the optimization procedure to fit the reduced operators $\hM$, $\hK$ to the reduced data, satisfying the \cref{eq:red_system}. The continuous-in-time displacement field is replaced by the reduced snapshots, defined as $\widehat{\Mat{Q}} = 
\begin{pmatrix}
\Mat{Q}_B\\
\widehat{\Mat{Q}}_I,
\end{pmatrix}$, where $\widehat{\Mat{Q}}_I \in \mathbb{R}^{r \times n_k}$. The corresponding details are covered in \cref{subsec:novel_optimization}.


\subsubsection{Reduction basis for the interior subsystem} \label{subsubsec:interior_basis}

By analogy with the Craig-Bampton method, the interior subsystem reduction basis $\bV_I$ should represent the dynamics of the interior subsystem. Originally, the truncated set of modes of the interior subsystem is used, but it requires access to the submatrices $\Mat{M}_{II}$ and $\Mat{K}_{II}$ that are not available in this work. Instead, we have to calculate the reduction basis from the snapshot data. 

Since we are free to define any admissible simulation conditions, we propose to create an additional snapshot set that would characterize the interior subsystem's response.
For this purpose, the original system has to be simulated under the same external load $\Vec{f} (t)$ that was used to produce the system response \eqref{mat:snap_Q}, but with different boundary conditions, namely \textit{with fixed boundary nodes}. This simulation corresponds to the solution of the following system:

\begin{equation} \label{eq:int_subsys}
\Mat{M}_{II}\ddot{\Vec{q}}^{(1)}_{I} (t) + \Mat{K}_{II}\Vec{q}^{(1)}_I (t) = \Vec{f}^{(1)}_I (t).
\end{equation}

We denote the state and force vector with the superscript $^{(1)}$ to distinguish these data from the snapshots \eqref{mat:snap_Q}. The system response is organized in the following snapshot matrix:

\begin{equation} \label{mat:red_int_snap}
    \bQ^{(1)} =
\begin{pmatrix}
| &  \dots & | \\
\Vec{q}^{(1)}_I(t_1) & \dots & \Vec{q}^{(1)}_I(t_N) \\
| & \dots & |
\end{pmatrix}.
\end{equation}.

The target basis $\Mat{V}_I$ is obtained by taking the first $r$ left singular vectors of the singular value decomposition of the $\bQ^{(1)}$.

Moreover, already using the current data, we are able to identify the reduced interior subsystem, applying the force-informed operator inference methodology, described in subsection \ref{subsec:fi_opinf}. 
 \begin{equation} \label{eq:finf_opt_I}
    \hM_{II}^{(1)}, \hK_{II}^{(1)} = \text{arg} \min_{\substack{\hM_{II} \succ 0 \\ \hK_{II} \succ 0}} \| \hM_{II} \ddot{\widehat{\textbf{Q}}}^{(1)}_I  + \hK_{II} \widehat{\textbf{Q}}^{(1)}_I - \widehat{\textbf{F}}^{(1)}_I \| _{_F}^2 .
 \end{equation}
Although the resulting reduced operators cannot be used as a final form of $\hM_{II}$ and $\hK_{II}$, because the influence of the boundary subsystem was not taken into account, they are very important to assure the compatibility of the overall ROM, which is explained in more details in \cref{subsec:novel_optimization}.

Next, we will discuss three different approaches to obtain the coupling matrix $\Mat{\Phi}_{IB}$, which is an important component for the accurate response of the ROM under contact conditions.

\subsubsection{Approximation of the coupling matrix} \label{subsubsec:coupling_mat}
\begin{figure}[H]
    \centering
    \includegraphics[width=0.5\linewidth]{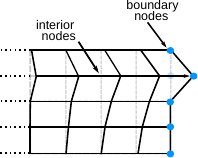}
    \caption{Static mode.}
    \label{fig:unit_displacements}
\end{figure}

In the following, we propose several ways to approximate the matrix $\boldsymbol{\Phi}_{IB}$ in a non-intrusive manner.
First of all, we consider the expression for the interior displacement vector given in \eqref{eq:red_q_relation}.
Since the ROM is not constructed yet, we can only approximate $\widehat{\Vec{q}}_I$ with the projected data from an interior subsystem simulation, as it is shown in \cref{subsubsec:interior_basis} by $\bq^{(1)}_I$. Using the snapshot matrix notation we replace the relation \eqref{eq:red_q_relation} with the following assumption:
\begin{equation} \label{eq:coupling_edited}
    \bQ_{I} = \boldsymbol{\Phi}_{IB} \bQ_{B} + \bQ_{I}^{(1)}
\end{equation} Thus the least-squares problem to obtain the coupling matrix can be formulated:

\begin{equation} \label{eq:lsq_z_full}
    \bs{\Phi}_{IB} = \text{arg} \min_{\Phi_{IB}} \| \bQ_{I} - \boldsymbol{\Phi}_{IB} \bQ_{B} + \bQ_{I}^{(1)} \| _{_F}^2 .
 \end{equation}
The disadvantage of this approach is the large dimension of the least-squares problem. To overcome this difficulty, we propose to reduce the number of unknowns by projecting the equation \eqref{eq:coupling_edited} onto the low-dimensional subspace spanned by the POD modes of the snapshot set $\bQ_I$, denoted by $\bV_I^{(2)}$. In this case, we can likewise identify the reduced form of the coupling matrix via the following least-squares problem:

\begin{equation} \label{eq:lsq_redz}
    \widetilde{\boldsymbol{\Phi}}_{IB} = \text{arg} \min_{\widetilde{\boldsymbol{\Phi}}_{IB}} \| \bV_I^{(2)}{}^T \bQ_{I}^{(2)} - \widetilde{\boldsymbol{\Phi}}_{IB} \bQ_{B}^{(2)} + \bV_I^{(2)}{}^T \bQ_{I}^{(1)} \| _{_F}^2 ,
 \end{equation} while the full-dimensional coupling matrix is then recovered as:

\begin{equation} \label{eq:z_full_from_red}
    \boldsymbol{\Phi}_{IB} = \bV^{(2)}_I \widetilde{\boldsymbol{\Phi}}_{IB} .
\end{equation} The additional projection with $\Mat{V}_I^{(2)}$ allows us to save memory and time resources for obtaining the coupling matrix. In \cref{sec:5} we compare the results for both types of the LSQ approximation of $\boldsymbol{\Phi}_{IB}$. 
Since the coupling matrix plays a very important role for the accurate approximation of the contact forces, we propose an additional straightforward method to compute the exact $\boldsymbol{\Phi}_{IB}$. For this purpose, we will use its original interpretation as the constrained modes of the system \cite{Craig1968}, \cite{subAllet2020}.
In other words, the corresponding columns of the coupling matrix are the interior system response when one boundary freedom is subjected to a unit displacement, while other boundary freedoms are fixed, see \cref{fig:unit_displacements}.

In this case, the identification of the coupling term is exact. If the assumption about the small size of the contact region holds, the calculation of the coupling term through consecutive static simulations is not only highly accurate but also very efficient.
However, as the number of contact nodes increases, the effectiveness of this method is decreasing drastically. 

\subsubsection{Optimization problem for system matrices} \label{subsec:novel_optimization} 
Once the interior reduction basis and the coupling term are determined, the actual inference of the reduced-order system can be completed. The global reduction basis can be then calculated as \eqref{mat:V}. 

The last missing link for the optimization problem analogous to \eqref{eq:finf_opt} is the low-dimensional data representation. Note, that the matrix $\Mat{V}$ is not orthogonal. As a consequence, it is not correct to calculate the low-dimensional data representation by multiplying the given snapshots $\Mat{Q}$ with the transpose of the global reduction basis, i.e.,  $\widehat{\Mat{Q}} \neq \Mat{V}^T \Mat{Q}$. Instead, it is necessary to consider the direct relation $\Mat{Q} = \Mat{V} \widehat{\Mat{Q}}$. The boundary degrees of freedom remain in their full dimension, whereas the reduced interior degrees of freedom could be represented using the subspace spanned by $\Mat{V}^{(1)}$, for example, calculated as $\Mat{V}^{(1) T} \Mat{Q}_I^{(1)}$, as it is summarized in \eqref{mat:final_red_data}.

The reduced force snapshot matrix, in turn, can be simply obtained as $\widehat{\Mat{F}} = \Mat{V}^T \Mat{F}$, which results in the corresponding expressions for its boundary and interior regions, as it is shown in \eqref{mat:final_red_data}. 
 
\begin{equation}
\label{mat:final_red_data}
\widehat{\textbf{Q}} =
\begin{pmatrix}
\textbf{Q}_B \\
\widehat{\textbf{Q}}_I{}^{(1)}
\end{pmatrix},  \;
\widehat{\textbf{F}} = \begin{pmatrix}
\textbf{F}_B + \bs{\Phi}_{IB} \bF_I^{(1)}  \\
\widehat{\textbf{F}}_I
\end{pmatrix} .
\end{equation}Note, that it holds $\Mat{F}_I = \Mat{F}_I^{(1)}$ by definition since we propose using the same external force for both simulations of \eqref{contactfreeODE} and \eqref{eq:int_subsys}.  The second derivatives of the state vector are calculated via the finite-difference method from the reduced snapshots $\widehat{\Mat{Q}}$.

Finally, the reduced operators can be inferred by solving an optimization problem. We include the SPD constraints, which are of significant importance also for the latter application to the contact mechanics problems. To have consistency with the reduced interior part of the system and the subspace spanned by $\Mat{V}_I^{(1)}$, we enforce the equality of the final reduced interior submatrix $\hM_{II}$ and the identified interior submatrix $\hM_{II}^{(1)}$ that was calculated in the subsection \ref{subsubsec:interior_basis}. The same holds for the stiffness matrix $\hK_{II}$. This leads us to the final optimization problem:

 \begin{equation} \label{eq:finf_opt_final}
    \hM, \hK = \text{arg} \min_{\substack{\hM \succ 0, \\ \hK \succ 0 \\ \hM_{II} = \hM_{II}^{(1)} \\ \hK_{II} = \hK_{II}^{(1)} }} \| \hM \ddot{\widehat{\textbf{Q}}}  + \hK \widehat{\textbf{Q}} - \widehat{\textbf{F}} \| _{_F}^2 .
 \end{equation} In the end, we are able to find the ROM using the data from the simulations with no contact constraints, preserving the SPD properties and the substructured form according to the original system operators. The resulting reduced model can be efficiently utilized not only for similar loading and boundary conditions covered in the training data but also for much more complex scenarios such as dynamic contact. In the next section, we will show the numerical results for two examples that prove the latter fact.

 \begin{algorithm}[H] 
\caption{Constrained substructured operator inference algorithm}\label{alg:final_inference}
\SetAlgoLined
\KwIn{$\bQ_B, \bQ_I \; \text{as in } \eqref{mat:snap_Q}, \bF_B, \bF_I \; \text{as in } \eqref{mat:snap_F}, \; \bQ^{(1)} \text{as in } \eqref{mat:red_int_snap}$}
\vspace{0.2cm}

\KwOut{$\widehat{\textbf{M}}, \widehat{\textbf{K}}$}
\vspace{0.2cm}

\nl The interior reduction basis and auxiliary interior sub-matrices $\textbf{V}_I , \; \hM^{(1)}_{II}, \; \hK^{(1)}_{II} \leftarrow$ compute using ($\bQ^{(1)} , \bF_I$), as described in \cref{subsubsec:interior_basis}
\vspace{0.2cm}

\nl The coupling matrix \\ $\bs{\Phi}_{IB} \leftarrow$ compute using ($\bQ^{(1)}, \; \bQ_B, \; \bQ_I$), as described in \cref{subsubsec:coupling_mat}.

\nl  The global reduction basis $\textbf{V} \leftarrow \begin{pmatrix}
			\textbf{I}_m & 0 \\
			\bs{\Phi_{IB}}   & \textbf{V}_r
\end{pmatrix}$
\vspace{0.2cm}

\nl  Reduced data $\widehat{\textbf{F}} \leftarrow \begin{pmatrix}
\textbf{F}_B + \bs{\Phi}_{IB} \bF_I^{(1)}  \\
\widehat{\textbf{F}}_I
\end{pmatrix} $, $\widehat{\textbf{Q}} \leftarrow
\begin{pmatrix}
\textbf{Q}_B^{(2)}  \\
\widehat{\textbf{Q}}_I{}^{(1)}
\end{pmatrix}$
\vspace{0.2cm}

\nl  The data matrix $\boldsymbol{\mathcal{D}} \leftarrow \begin{pmatrix}
     \ddot{\widehat{\textbf{Q}}}^T, &  \widehat{\textbf{Q}}^T
\end{pmatrix} $
\vspace{0.2cm}

\nl  The right-hand-side matrix $\boldsymbol{\mathcal{R}} \leftarrow 
     \widehat{\textbf{F}}^T$
\vspace{0.2cm}

\nl  Infer the reduced system: 
$ \widehat{\textbf{M}}, \widehat{\textbf{K}} \leftarrow \text{minimize} \lVert \boldsymbol{\mathcal{D}} \begin{pmatrix}
			\widehat{\textbf{M}}, & 			\widehat{\textbf{K}}
			\end{pmatrix} \,^T - \boldsymbol{\mathcal{R}} \rVert_{_F}^{^2}$, \\ \phantom{Infer the reduced system} s.t. $\widehat{\textbf{M}}  \succ 0 , \quad \widehat{\textbf{K}} \succ 0, \quad \hM_{II} = \hM_{II}^{(1)},  \quad \hK_{II} = \hK_{II}^{(1)}$

\vspace{0.2cm}

\end{algorithm}

\subsection{Solving the adjoint system} \label{subsec: adjoint}

After inferring the reduced primal system with Alg.~\ref{alg:final_inference}, the contact constraints are appended. Due to contact treatment by substructuring, the contact constraints remain unaffected by the reduction procedure and preserve their structure. In particular, they can be described as follows:
\begin{align} \label{ineq_red_1}
    \Mat{C}_B\widehat{\bq}_B + \Vec{b} \geq \Vec{0}, \quad \Vec{\lambda} \geq \Vec{0}, \quad \Vec{\lambda}^T (\Mat{C}_B\widehat{\bq}_B + \Vec{b})  = \Vec{0},
\end{align}
where $\Mat{C}_B$ corresponds to the boundary block-matrix of $\Mat{C}.$
The variable $\widehat{\bq}_B$ in \eqref{ineq_red_1} corresponds to the boundary displacement nodes of the reduced variable $\widehat{\bq},$ which can be identified due to the Craig-Bampton partitioning and represents a direct approximation for the boundary displacements of the full order model.

In the following, we consider equidistant time grid $\{t_i: i = 0,\ldots, k\}$ with $t_0 = 0$ and $t_{k} = T.$ We denote $\widehat{\Vec{q}}(t_i) = \widehat{\Vec{q}}_i$ and $\Vec{\lambda}(t_i) = \Vec{\lambda}_i.$ The time step is denoted as $h.$

Applying now the primal-dual decoupling approach from \cref{subsec:primaldual} to the reduced primal system, for $i$th time step with $i > 1,$ we obtain the following LCP

\begin{equation} \label{eq:red_KKTAB}
	\begin{array}{rcc}
	\widehat{\Vec{B}}_i + \widehat{\Mat{A}} \Vec{\lambda}_i & \geq & 0, \\
	\Vec{\lambda}_i & \geq & 0, \\
	\Vec{\lambda}_i^T (\widehat{\Vec{B}}_i + \widehat{\Mat{A}} \Vec{\lambda}_i) & = & 0.
	\end{array}
\end{equation}

where

\begin{eqnarray} \label{mat:redLCP}
\widehat{\Mat{A}} &:=& h^2 \Mat{C}_{B}(\widehat{\Mat{M}}+h^2 \widehat{\Mat{K}})^{-1}\Mat{C}_B^T,\\
\widehat{\Mat{B}}_i &:=& \Mat{C}_B (\widehat{\Mat{M}}+h^2 \widehat{\Mat{K}})^{-1}(h^2 \widehat{\Vec{f}}_i + 
        2 \widehat{\Mat{M}} \widehat{\Vec{q}}_{i-1} - \widehat{\Mat{M}} \widehat{\Vec{q}}_{i-2}) + \Vec{b}
\end{eqnarray}

When the Lagrange multiplier is known, the reduced displacement vector, at the corresponding time step can be computed by:

\begin{equation}\label{euler_ntn_red}
\begin{aligned}
\widehat{\Vec{q}}_i = \big(\widehat{\Mat{M}}+h^2\widehat{\Mat{K}}\big)^{-1}
(h^2 \widehat{\Vec{f}}_i + h^2\Mat{C}_B^T\widehat{\Vec{\lambda}}_i +
 2\widehat{\Mat{M}}\widehat{\Vec{q}}_{i-1}- \widehat{\Mat{M}}\widehat{\Vec{q}}_{i-2}).
 \end{aligned}
\end{equation}

\begin{algorithm}[H] 
    \caption{\label{alg:contactscheme} Solution algorithm of reduced node-to-node contact problems}
    \SetAlgoLined
    \DontPrintSemicolon
    \KwIn{$\bQ_B, \bQ_I \; \text{as in } \eqref{mat:snap_Q}, \bF_B, \bF_I \; \text{as in } \eqref{mat:snap_F}, \; \bQ^{(1)} \text{as in } \eqref{mat:red_int_snap} , \Mat{C}_B, \Vec{b} \text{ as in } \eqref{ineq_red_1} $} 
    \vspace{0.2cm}
	\KwOut{$\Vec{q}_{i}$ and $\Vec{\lambda}_{i}$  for $i \in \{2,\ldots,n_T\}$}
    \vspace{0.2cm}
    \nl Infer the matrices $\widehat{\Mat{M}}, \widehat{\Mat{K}}$ by means of Alg.~\ref{alg:final_inference}.\;
    \nl Append the contact constraints defined by $\Mat{C}_B, \Vec{b}.$\;
    \nl Compute the LCP matrix $\widehat{\Mat{A}}$ in \eqref{mat:redLCP}\;
    \vspace{0.1cm}
	\nl Initialization of $\widehat{\Vec{q}}_i$ for $i = 0,1.$\;
	 \For{$i=2,\dots, n_T$}{
        \nl Compute $\widehat{\Mat{B}}_i$ by \eqref{mat:redLCP}\;
        \vspace{0.1cm}
        \nl Compute $\widehat{\Vec{\lambda}}_{i}$ as a solution of the LCP \eqref{eq:red_KKTAB} by Lemke's method.\;
        \vspace{0.1cm}
    	\nl Compute the reduced displacement vector $\widehat{\Vec{q}}_{i}$ by \eqref{euler_ntn_red}.}
     \vspace{0.1cm}
	\nl Compute $\Vec{q}_{i} = \Mat{V}\widehat{\Vec{q}}_{i},\ i = 2,\ldots n_T.$
	 	
\end{algorithm}
\vspace{0.2cm}

The solving procedure for the reduced contact problem is presented in Alg.~\ref{alg:contactscheme}. Here, the LCP \eqref{eq:red_KKTAB} is solved for each time step by Lemke's method \cite{LemkeAlg} providing the Lagrange multipliers $\widehat{\Vec{\lambda}}_i.$ Since the SPD properties of the reduced matrices are preserved, the properties are transferred to the LCP matrix $\widehat{\Mat{A}}$ in \eqref{mat:redLCP}. According to \cite{CotPetal2009}, the LCP \eqref{mat:redLCP} has a unique solution for all $\widehat{\Mat{B}}_i$, if the LCP matrix $\widehat{\Mat{A}}$ is SPD. Therefore, the LCP in Alg.~\ref{alg:contactscheme} will converge.

Unlike the reduced displacement vector, the Lagrange multiplier $\widehat{\Vec{\lambda}}_i$ does not require an initial value and thus it possesses no memory. Even more, it preserves its dimension and serves as a direct approximation of the Lagrange multipliers computed by the full-order model. Moreover, since the shape of the contact region depends on the acting force, this shape has to be computed also in the reduced model. It follows that the performance of our algorithm depends on the ratio of contact nodes and interior volume nodes, because only the latter may be reduced. Further improvements can be achieved by introducing strategies to keep the maximal number of contact nodes small. One possibility would be to provide in the offline phase two reduced models: one with a complete set of possible contact nodes and one with only a small subset of contact nodes. In the online phase, we could choose the appropriate reduced model depending on the evolution of the contact area.

\section{Numerical Experiments} \label{sec:5}
In the following, we present the performance of the methodology described in \cref{sec:4}.  For this purpose, two numerical examples have been used: a structural beam with a rigid plane underneath it, and a connecting rod with an obstacle restricting the movement of its crank pin. The CAD geometries of both models were created using the
Simcenter 3D software. All the algorithms defined in Section 5.2 are implemented in a custom Python environment. Time integration is performed with the implicit Euler scheme \cite{Bel83}. The test simulations with active contact conditions are performed using Lemke's algorithm \cite{CotPetal2009} for solving LCP problems from an open-source Python library \cite{LemkeAlg}. The semidefinite programming algorithms from the CVXPY library \cite{diamond2016cvxpy} are used to solve the constrained optimization problems \eqref{eq:finf_opt_I} and \eqref{eq:finf_opt_final}.

In the following study, the variables of interest are the displacement vector and the Lagrange multipliers. The constraint matrix $\Mat{C}$ is formulated as a dimensionless Boolean matrix. Consequently, the computed Lagrange multipliers correspond to the contact force values. 

To measure the performance of the models, we used the plots of the trajectories as well as the relative errors calculated as: \begin{align}\label{eq:relerror}
    \bs{\varepsilon}_q = \frac{ \| \Vec{q} (t) - \hq(t) \|_F^2 }{\max \| \Vec{q} (t) \|_F^2}, & \quad \bs{\varepsilon}_{\lambda} = \frac{ \| \Vec{\lambda} (t) - \widehat{\Vec{\lambda}} (t) \|_F^2 }{\max \| \Vec{\lambda} (t) \|_F^2}.
\end{align} At every time step, the error is divided by the maximum norm of the original trajectory to avoid drastically increasing the error near the zero level.

\subsection{Console structural beam with an obstacle} \label{beamexample}
\begin{figure}[H]
     \centering
     \begin{subfigure}[b]{0.4\textwidth}
         \centering
         \includegraphics[width=\textwidth]{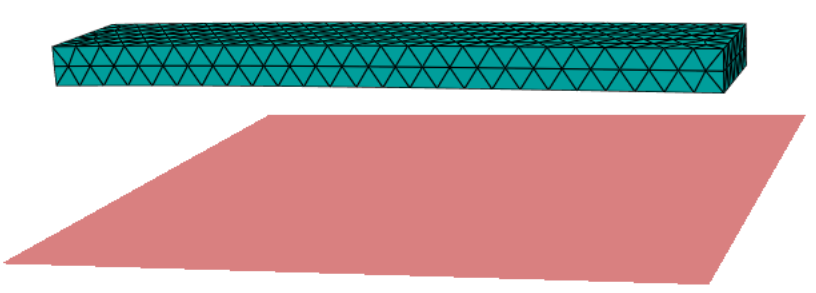}
         \caption{Undeformed shape.}
         \label{fig:beam_a}
     \end{subfigure}
     \hfill
     \begin{subfigure}[b]{0.4\textwidth}
         \centering
         \includegraphics[width=\textwidth]{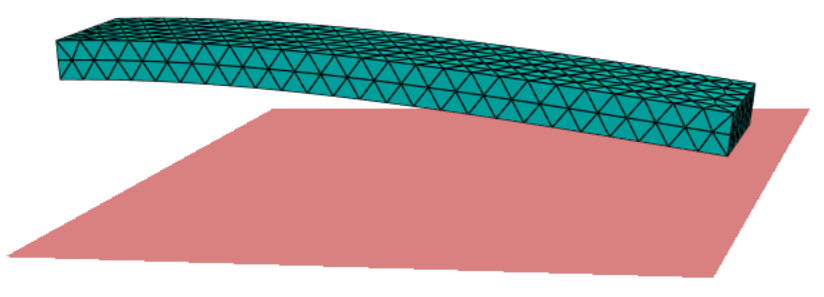}
         \caption{Deformed shape.}
         \label{fig:beam_b}
     \end{subfigure}
     \hfill
        \caption{FEM model of a console beam with a rigid plane obstacle.}
        \label{fig:beam}
\end{figure}

Our first numerical study concerns a standard 3D model of a linear elastic structural beam discretized by a finite element grid with 1599 degrees of freedom. The corresponding material properties are: Young's modulus $E$ = 210 GPa, Poisson's ratio $\nu$ = 0.3, and density $\rho$ = 7860 $\mathrm{kg} / \mathrm{m}^3$. The beam is fixed on the left side, and a rigid plane is placed one meter below the beam. The force is applied along the line at the free end of the beam. Fig.~\ref{fig:beam_a} shows the FEM model of the undeformed beam, and Fig.~\ref{fig:beam_b} depicts the corresponding 3D shape after loading, including the contact with the obstacle.  There are six nodes where the beam touches the plane as soon as the deflection enables contact. These nodes are selected as boundary degrees of freedom $\Vec{q}_B$. All other nodes represent interior freedoms $\Vec{q}_I$. 

\begin{figure}[H]
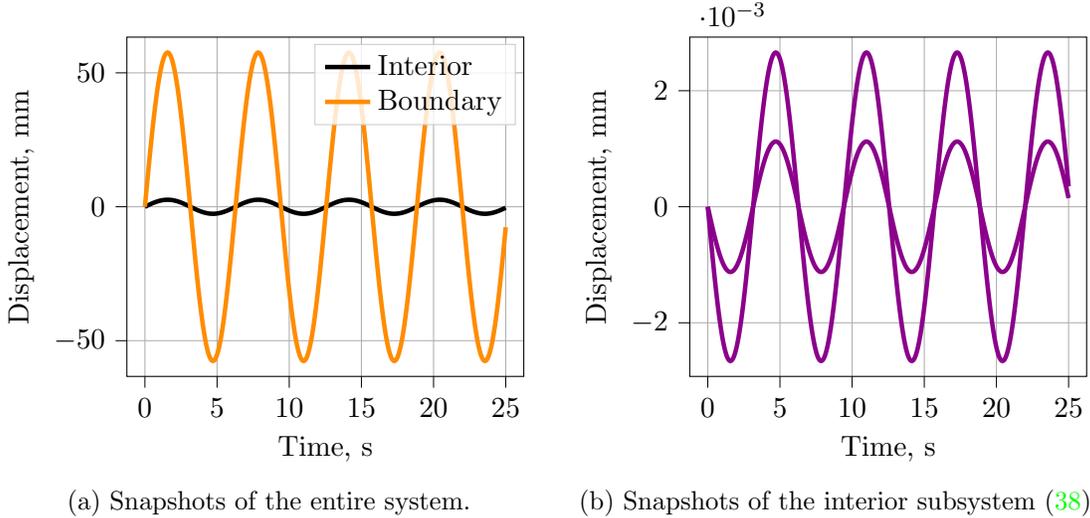

     \centering
     \begin{subfigure}[b]{0.49\textwidth}
         \centering
         \input{plots/training_data_snapshots}
         \caption{Snapshots of the entire system.}
         \label{fig:train_a}
     \end{subfigure}
     \hfill
     \begin{subfigure}[b]{0.49\textwidth}
         \centering
         \input{plots/training_data_interiorQ_snapshots}
         \caption{Snapshots of the interior subsystem \eqref{eq:int_subsys}.}
         \label{fig:train_b}
     \end{subfigure}
     \hfill
        \caption{Training snapshots for the console beam example.}
        \label{fig:train}
\end{figure}

To construct non-intrusive reduced-order models we used simulations with the harmonic external force with a frequency $0.16$ Hz and an amplitude of 3 kN. The training data are represented in Fig.~\ref{fig:train}, where Fig.~\ref{fig:train_a} shows selected snapshots used for the coupling matrix approximation (\cref{subsubsec:coupling_mat}) and the final inference problem (Alg.~\ref{alg:final_inference}), while Fig.~\ref{fig:train_b} shows the snapshots used for the interior reduction basis calculation (\cref{subsubsec:interior_basis}).

To get an estimate of the reduced order of the model, we analyze the singular value decay of the given interior snapshots displayed in Fig.~\ref{fig:train}. According to a rapid decay of singular values $\sigma$ in Fig.~\ref{fig:svd_beam}, we have chosen the dimension $r = 2$ for the reduction basis of the interior part $\Mat{V}_I^{(1)}$ and for the reduction basis for the coupling approximation $\Mat{V}_I^{(2)}$.  \begin{figure}[H]
     \centering
         \centering
\begin{tikzpicture}

\definecolor{darkgray176}{RGB}{176,176,176}
\definecolor{darkmagenta}{RGB}{139,0,139}
\definecolor{lightgray204}{RGB}{204,204,204}

\begin{axis}[
width=2.7in,
height=2.4in,
legend cell align={left},
legend style={at={(axis cs:2.0,1e-8)},anchor=south west, cells={align=left}, fill opacity=0.8, draw opacity=1, text opacity=1, draw=lightgray204},
log basis y={10},
tick align=outside,
tick pos=left,
x grid style={darkgray176},
xmajorgrids,
xmin=-1.45, xmax=30.45,
xtick style={color=black},
y grid style={darkgray176},
ymajorgrids,
ylabel={ \small Normalized $\sigma$},
xlabel={ \small $r$},
ylabel style={yshift=-0.1cm},
ymin=1.58096622064817e-17, ymax=6.31031886771027,
ymode=log,
ytick style={color=black},
ytick={1e-19,1e-17,1e-15,1e-13,1e-11,1e-09,1e-07,1e-05,0.001,0.1,10,1000},
yticklabels={
  \(\displaystyle {10^{-19}}\),
  \(\displaystyle {10^{-17}}\),
  \(\displaystyle {10^{-15}}\),
  \(\displaystyle {10^{-13}}\),
  \(\displaystyle {10^{-11}}\),
  \(\displaystyle {10^{-9}}\),
  \(\displaystyle {10^{-7}}\),
  \(\displaystyle {10^{-5}}\),
  \(\displaystyle {10^{-3}}\),
  \(\displaystyle {10^{-1}}\),
  \(\displaystyle {10^{1}}\),
  \(\displaystyle {10^{3}}\)
}
]
\addplot [semithick, darkmagenta, mark=*, mark size=3, mark options={solid}, only marks]
table {%
0 1
1 2.14381119755648e-06
2 1.66858967378909e-10
3 3.2783923921379e-13
4 2.60799419868915e-14
5 1.99963367023467e-15
6 1.4759513017458e-15
7 7.02208660548218e-16
8 1.73408668593413e-16
9 1.70759866675055e-16
10 1.09214970361769e-16
11 1.04060040374281e-16
12 9.97640097136873e-17
13 9.97640097136873e-17
14 9.97640097136873e-17
15 9.97640097136873e-17
16 9.97640097136873e-17
17 9.97640097136873e-17
18 9.97640097136873e-17
19 9.97640097136873e-17
20 9.97640097136873e-17
21 9.97640097136873e-17
22 9.97640097136873e-17
23 9.97640097136873e-17
24 9.97640097136873e-17
25 9.97640097136873e-17
26 9.97640097136873e-17
27 9.97640097136873e-17
28 9.97640097136873e-17
29 9.97640097136873e-17
};
\addlegendentry{\small Interior subsystem sim.}
\addplot [semithick, black, mark=*, mark size=3, mark options={solid}, only marks]
table {%
0 1
1 2.18603631275728e-06
2 1.47333644428407e-10
3 4.96498994703349e-13
4 7.17169732515814e-15
5 1.40025198718262e-15
6 6.42452221855788e-16
7 2.67911576003116e-16
8 1.71552677724878e-16
9 1.30856034969627e-16
10 9.99200635313693e-17
11 9.99200635313693e-17
12 9.99200635313693e-17
13 9.99200635313693e-17
14 9.99200635313693e-17
15 9.99200635313693e-17
16 9.99200635313693e-17
17 9.99200635313693e-17
18 9.99200635313693e-17
19 9.99200635313693e-17
20 9.99200635313693e-17
21 9.99200635313693e-17
22 9.99200635313693e-17
23 9.99200635313693e-17
24 9.99200635313693e-17
25 9.99200635313693e-17
26 9.99200635313693e-17
27 9.99200635313693e-17
28 9.99200635313693e-17
29 9.99200635313693e-17
};
\addlegendentry{ \hspace{0.2cm} \\ \small Interior part of the \\ \small full-system sim.}
\end{axis}

\end{tikzpicture}
         \vspace*{-0.4cm}
         \caption{Console beam example: SVD decay of the snapshot set that corresponds to the interior part of the full-order system simulation, and the interior subsystem snapshot set of the full-order system simulation with fixed boundary nodes.}
         \label{fig:svd_beam}
\end{figure}
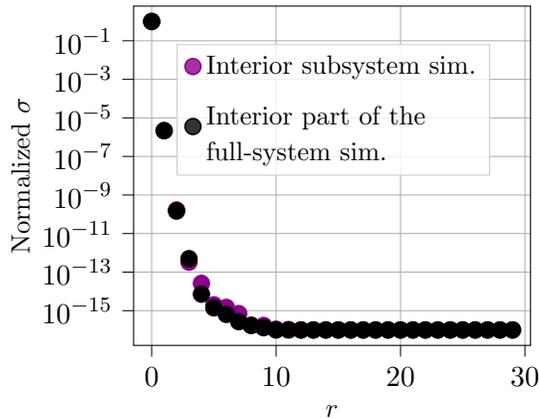 The figures below show the response of three ROMs of order $r = 2$ derived with different coupling matrix determination strategies that have been proposed in \cref{subsubsec:coupling_mat}. The models were subjected to a harmonic force, and in addition, the rigid surface at the bottom of the beam forms an obstacle, creating contact stresses, when the beam reaches the surface. The load frequency has been set higher than the frequency used to generate snapshots, namely to $0,32$ Hz. The system response is calculated by solving the full-dimensional contact problem \eqref{eq:KKT_transient1} and \eqref{eq:KKT_transient2}, and the reduced-dimensional LCP \eqref{eq:red_KKTAB}. 

\begin{figure}[H]
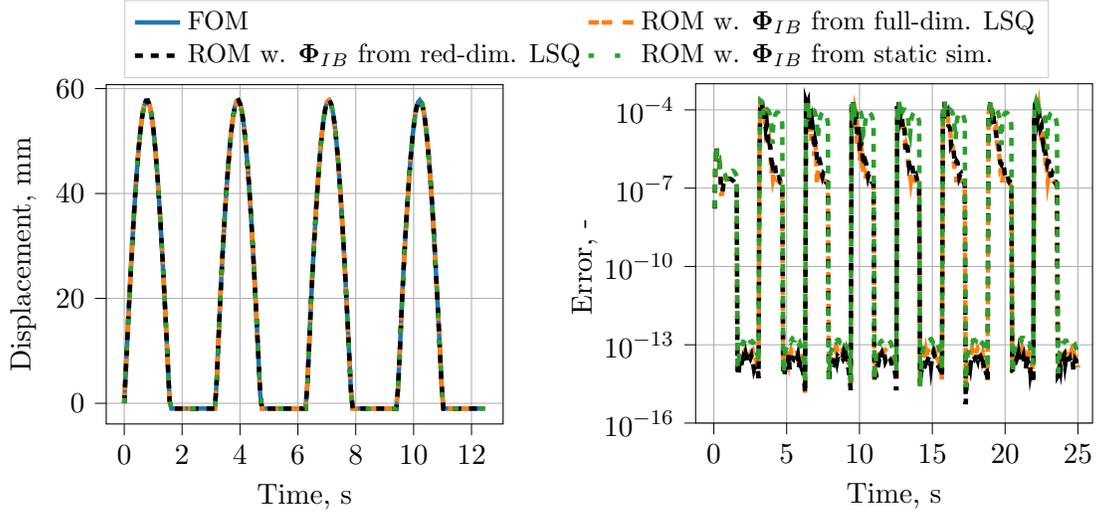

     \centering
     \begin{subfigure}[b]{0.47\textwidth}
         \centering
         \input{plots/UB_trajectories}
         \vspace*{-0.3cm}
         \caption{Boundary displacement trajectory $\Vec{q}_B$.}
         \label{fig:bound_a}
     \end{subfigure}
     \hfill
     \begin{subfigure}[b]{0.52\textwidth}
         \centering
         \input{plots/UB_errors}
         \vspace*{-0.3cm}
         \caption{Relative errors for the boundary displacements.}
         \label{fig:bound_b}
     \end{subfigure}
        \caption{Console beam model: Comparison of the contact simulation results $\Vec{q}_B$ of the FOM and substructured OpInf ROMs with different coupling term approximations.}
        \label{fig:bound}
\end{figure}

\begin{figure}[H]
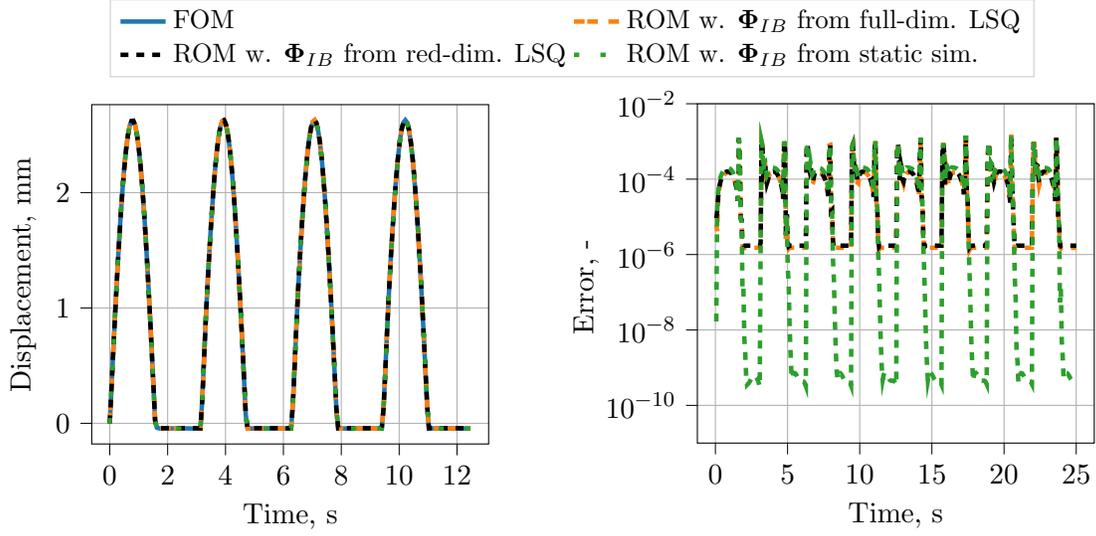

     \centering
     \begin{subfigure}[b]{0.47\textwidth}
         \centering
         \input{plots/UI_trajectories}
         \vspace*{-0.4cm}
         \caption{Interior displacement trajectory $\Vec{q}_I$.}
         \label{fig:inter_a}
     \end{subfigure}
     \hfill
     \begin{subfigure}[b]{0.52\textwidth}
         \centering
         \input{plots/UI_errors}
         \vspace*{-0.4cm}
         \caption{Relative errors for the interior displacements.}
         \label{fig:inter_b}
     \end{subfigure}
        \caption{Console beam model: Comparison of the contact simulation results $\Vec{q}_I$ of the FOM and substructured OpInf ROMs with different coupling term approximations.}
        \label{fig:inter}
\end{figure}

Fig.~\ref{fig:bound_a} shows the displacement trajectory of one of the boundary nodes. All reduced model solution curves match the original trajectory, which holds for every boundary node as proved by the error curves in Fig.~\ref{fig:bound_b}. The relative errors are calculated by equation \eqref{eq:relerror}. Similar outcomes occur for the interior nodes of the model, see Fig.~\ref{fig:inter}. Although the maximum relative error for the interior displacements of all the inferred ROMs is larger than for the boundary displacements, it still does not exceed 1 \% , which provides sufficient accuracy for the approximation of the system displacements.


\begin{figure}[H]
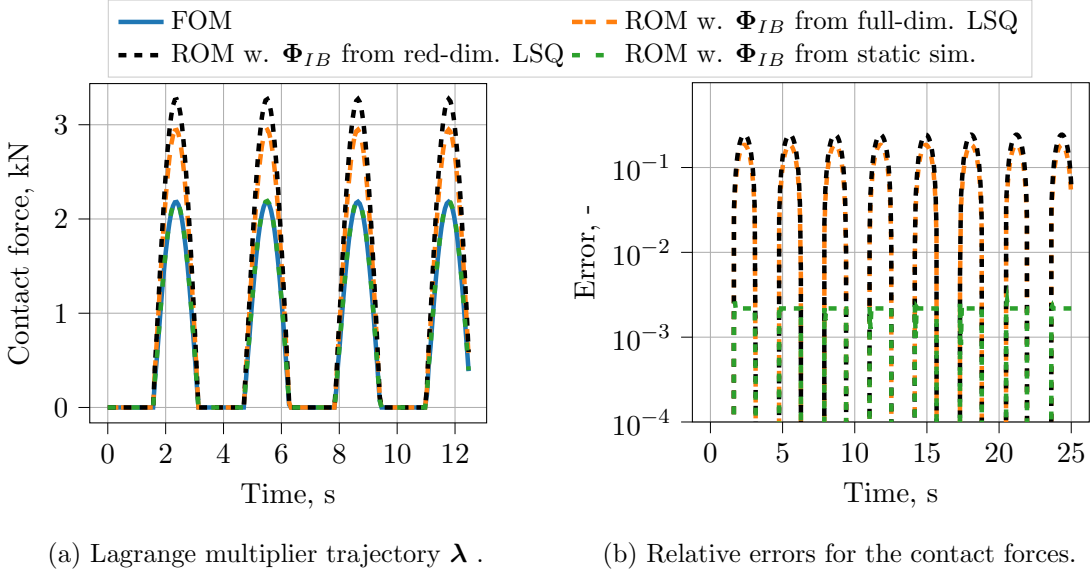

     \centering
     \begin{subfigure}[b]{0.47\textwidth}
         \centering
         \input{plots/LMB_trajectories}
         \vspace*{-0.4cm}
         \caption{Lagrange multiplier trajectory $\Vec{\lambda}$ .}
         \label{fig:lmb_a}
     \end{subfigure}
     \hfill
     \begin{subfigure}[b]{0.52\textwidth}
         \centering
         \input{plots/LMB_errors}
         \vspace*{-0.4cm}
         \caption{Relative errors for the contact forces.}
         \label{fig:lmb_b}
     \end{subfigure}
        \caption{Console beam model: Comparison of the contact simulation results $\Vec{\lambda}$ of the FOM and substructured OpInf ROMs with different coupling term approximations.}
        \label{fig:lmb}
\end{figure}

The most distinctive comparison corresponds to the numerical results for the Lagrange multipliers. As indicated in Fig.~\ref{fig:lmb}, the error is large for the ROMs that use the least-squares approximation of the coupling matrix. Although these reduced models correctly capture the qualitative behavior of the contact force, the numerical accuracy is not acceptable enough, and the contact force is predicted to be much higher than that of the original FOM. However, in the case when the coupling matrix is calculated "exactly", i.e., from the consecutive static simulations, the relative error drops to less than 1\%, and the resulting trajectory overlaps the contact force trajectory of the FOM, see Fig.~\ref{fig:lmb_a}.

Thus, the numerical results for the console beam example show the crucial importance of the correct coupling matrix approximation to represent the contact forces and corresponding contact pressures and stresses. Only the model with the exact coupling matrix has shown accurate contact force distribution.

However, all the three ROMs lead to a good approximation of the displacement field in the contact region. Moreover, the ROMs with least-squares approximations of the coupling matrix can still provide physically meaningful system response, including the qualitatively correct contact forces distribution, which could be used for the contact-search problems. For complex structures, it is often not known in advance if all the nodes in the defined contact region will actually become active, i.e., the contact nodes are just a guess. Hence, to define the contact region more precisely and at a low computational cost, one can use the ROM with LSQ approximation of the coupling matrix.

In conclusion, the proposed strategies can approximate the contact behavior of the console beam model, even though no information about the contact pressures was contained in the snapshots used for the non-intrusive ROM construction.

\subsection{Piston Rod with an obstacle}
\begin{figure}
     \centering
     \includegraphics{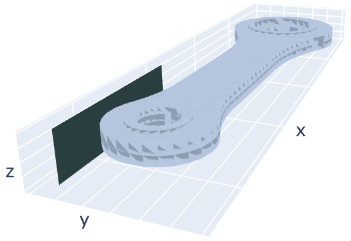}
     \vspace{0.1cm}
     \caption{A piston rod with a wall-obstacle.}
     \label{fig:pistonrod1}
 \end{figure}

In this section, we present the results of numerical experiments, which were conducted to analyze the behavior and performance of the piston rod mechanism. The piston rod, a crucial component in internal combustion engines, connects the piston to the crankshaft. Here we assume a linear Kirchhoff material law with Young's modulus $E$ = 210 GPa, Poisson's ratio $\nu$ = 0.3, and density $\rho$ = 7860 $\text{kg}/ \text{m}^3$. The piston rod contains two inner circles: a larger fixed circle and a smaller movable one. The larger circle is subjected to Dirichlet boundary conditions, while the smaller circle is free to move. An oscillating surface load is applied to the smaller circle, inducing movement in the corresponding part of the piston rod within the $x\text{-}y$ plane, see Fig.~\ref{fig:pistonrod1}. 
 \FloatBarrier
We first simulate the original system in a contact-free setting to test our novel approach from \Cref{sec:4}. Next, an additional obstacle in the form of a wall is introduced to address the extended contact problem and to compute the corresponding solution using the reduced model. The latter is then compared to the reference solution of the full-dimensional contact problem.

Similar to the experiment from Section~\ref{beamexample}, three solution variables of ROM and FOM are computed and compared: the interior displacements, the boundary displacements, and their counterparts, the Lagrange multipliers.

 \begin{figure}[H]
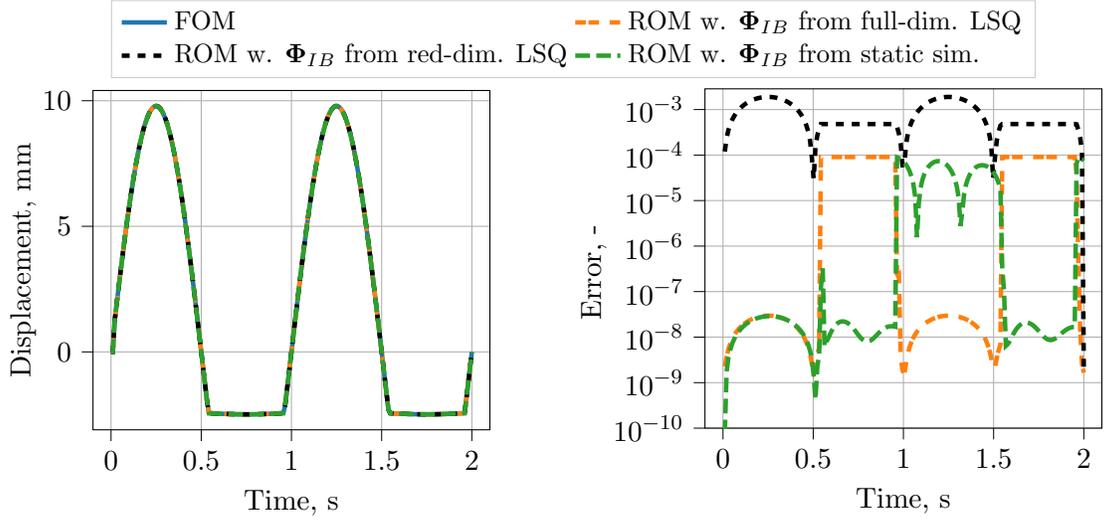

     \centering
     \begin{subfigure}[b]{0.47\textwidth}
         \centering
         \input{plots/pistonrod_slave_1}
         \vspace*{-0.5cm}
         \caption{Interior displacement trajectory $\Vec{q}_I$.}
         \label{fig:piston_slave_a}
     \end{subfigure}\hfill
     \begin{subfigure}[b]{0.52\textwidth}
         \centering
         \input{plots/pistonrod_error_slave_1}
         \caption{Relative errors for the interior displacements.}
         \label{fig:piston_slave_b}
     \end{subfigure}
        \caption{Piston rod model: Comparison of the contact simulation results $\Vec{q}_I$ of the FOM and substructured OpInf ROMs with different coupling term approximations.}
        \label{fig:piston_slave}
\end{figure}

 \begin{figure}[H]
     \centering
     \begin{subfigure}[b]{0.47\textwidth}
         \centering
         \input{plots/pistonrod_master_1}
         \vspace*{-0.5cm}
         \caption{Boundary displacement trajectory.}
         \label{fig:piston_master_a}
     \end{subfigure}\hfill
     \begin{subfigure}[b]{0.52\textwidth}
         \centering
         \input{plots/pistonrod_error_master_1}
         \caption{\small Relative errors for the boundary displacements.}
         \label{fig:piston_master_b}
     \end{subfigure}
        \caption{Piston rod model: Comparison of the contact simulation results $\Vec{q}_B$ of the FOM and substructured OpInf ROMs with different coupling term approximations.}
        \label{fig:piston_master}
\end{figure}

\begin{figure}[H]
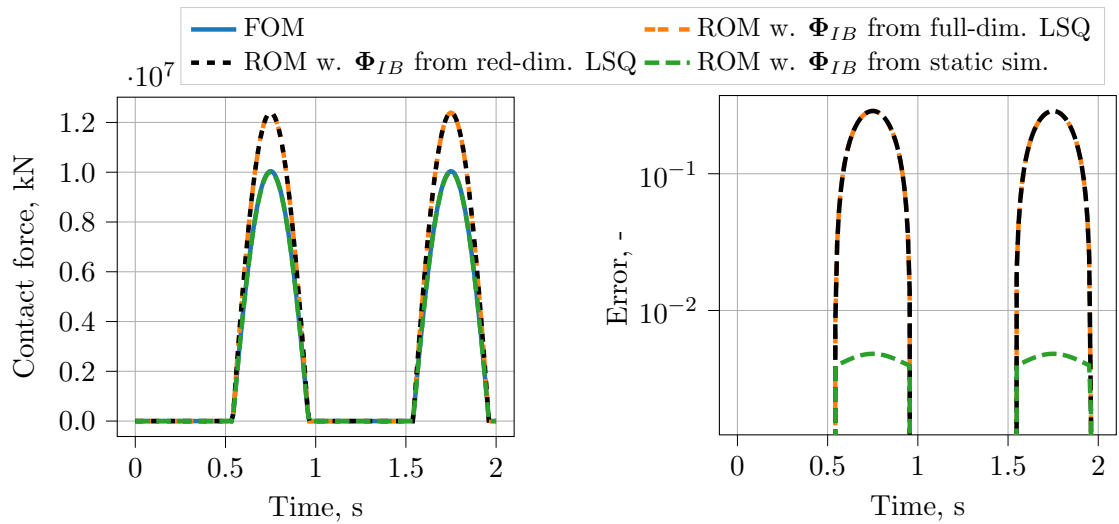

     \centering
     \begin{subfigure}[b]{0.49\textwidth}
         \centering
         \input{plots/pistonrod_lmb_1}
         \vspace*{-0.5cm}
         \caption{The Lagrange multipliers.}
         \label{fig:pistonrod_lmb_a}
     \end{subfigure}\hfill
     \begin{subfigure}[b]{0.49\textwidth}
         \centering
         \input{plots/pistonrod_error_lmb_1}
         \caption{Relative errors for the Lagrange multipliers.}
         \label{fig:pistonrod_lmb_b}
     \end{subfigure}
        \caption{Piston rod model: Comparison of the contact simulation results $\Vec{\lambda}$ of the FOM and substructured OpInf ROMs with different coupling term approximations.}
        \label{fig:pistonrod_lmb}
\end{figure}

Similar to the beam example, our novel ROM accurately predicts the contact behaviour of the piston rod. For all three approaches for the computation of the coupling matrix, the displacements for both the interior and the boundary nodes are well predicted. When observing the relative error trajectories, it is evident that the relative error of the boundary displacement nodes is even less than that of the interior displacement nodes, which agrees well with the fact that we preserve the contact shape and the dimension of the boundary displacements.

Regarding the results for the Lagrange multiplier, it is clear that the third approach of computing the coupling matrix achieves greater accuracy compared to the first two scenarios. However, in all three cases, the main two important aspects of contact pressure are recovered: The shape of the Lagrange multiplier trajectories and the timing of the active (i.e., strictly positive) and non-active (i.e., zero) phases.

In summary, we can again emphasize that the proposed strategies approximate the contact behavior of both numerical examples well, even though no information about the contact pressures was contained
in the snapshots used for the non-intrusive ROM construction.

\newpage

\section{Conclusions} \label{sec:6}
In this work, we have introduced an operator-inference-based method for developing reduced-order models for node-to-node contact problems utilizing contact-free finite element analysis data. The only required information about the contact condition is the maximum set of the potential contact nodes, called boundary nodes. In particular, our novel data-driven methodology utilizes substructuring, analogous to the classical Craig-Bampton approach. By distinguishing between the boundary and interior displacement nodes,  two different simulation scenarios are generated. The first simulation is performed with free-moving and the second one with fixed boundary nodes. Both sets of snapshots are used to infer the primal system and do not contain further information about the contact conditions.

Special attention is devoted to recovering the coupling matrix between the boundary and interior nodes, since it plays a critical role in the accurate representation of the Lagrange multipliers, i.e., contact pressures. We have proposed three different approaches for the coupling term identification, two of them are based on a special least-squares problem, either in full-dimensional or reduced-dimensional space. The third approach requires additional static simulation data when unit displacements are successively prescribed for each of the boundary nodes. The resulting interior static modes form the columns of the coupling matrix. Depending on the size of the contact area, the third approach can be more cumbersome than the first two, however, it consistently succeeds in providing accurate Lagrange multipliers.  All three approaches are similarly efficient in the case of a small contact area.

The final inference of the reduced matrices is performed using the operator inference procedure with additional linear matrix inequality constraints, which guarantees the symmetry and positive definiteness of the resulting operators.

After obtaining the reduced-order model, node-to-node contact constraints are appended. The inferred reduced-order model is solved with adjoint methods, i.e., with a primal-dual decoupling method, which allows switching from the reduced primal to the dual system, described by LCP. The Lagrange multipliers are computed by Lemke's method, a well-established pivoting method, which leads to the unique solution, due to the preserved system operator properties. Since the efficiency of our algorithms depends on the size of the LCP, which is solved in each time iteration, a small number of contact nodes is advantageous here as well.

The performance of the proposed algorithm is verified by means of two three-dimensional use cases: a structural beam and a steel piston rod. Further work is investigated for treating more complex contact scenarios such as node-to-segment contact problems, as well as considering friction in the contact condition.

\subsection*{Data availability statement}
The simulation data and the code used to generate the results presented in this paper are openly available at \url{https://gitlab.mpi-magdeburg.mpg.de/filanova/contact_opinf}. The numerical experiments were run on a Dell computer with 16 x 12th Gen Intel(R) Core(TM) i5-12600K, 31GB RAM, 64Bit, running Ubuntu version 20.04.4. Small numerical differences in the optimization results across the machines can be observed due to differences in floating-point arithmetic, CPU architecture, BLAS backend, etc.

\section*{Acknowledgements}
The authors acknowledge the support and computational resources provided by the Max Planck Institute for Dynamics of Complex
Technical Systems, Magdeburg and Siemens AG, Munich. This research was also supported by the Research Training Group "Mathematical Complexity
Reduction", which is a Graduiertenkolleg (DFG RTG 2297) funded by Deutsche Forschungsgemeinschaft (DFG), Germany.

\bibliographystyle{abbrv_1}
\bibliography{references}%

\end{document}